\documentclass[11pt]{article}

\advance\textheight -\headheight
\advance\textheight -\headsep
\oddsidemargin\paperwidth
\advance\oddsidemargin -\textwidth
\divide\oddsidemargin2 
\ifdim\oddsidemargin<.5truein \oddsidemargin.5truein \fi
\advance\oddsidemargin -1truein
\evensidemargin\oddsidemargin
\topmargin\paperheight \advance\topmargin -\textheight
\advance\topmargin -\headheight \advance\topmargin -\headsep
\divide\topmargin2
\ifdim\topmargin<.5truein \topmargin.5truein \fi
\advance\topmargin -1.4truein\relax

\usepackage[T2A]{fontenc}
\usepackage[utf8]{inputenc}
\usepackage[pdftex]{graphicx}
\usepackage{amsmath}
\usepackage{amsfonts}
\usepackage{comment}
\usepackage{pgf,tikz}
\usepackage{float}
\usepackage{makecell}
\usetikzlibrary{positioning,arrows}
\usetikzlibrary[patterns]

\usepackage{caption}
\usepackage{subcaption}
\usepackage{array}

\newcommand{\HH}{{\cal H}}
\newcommand{\CC}{\mathbb C}
\newcommand{\PP}{\mathbb P}
\newcommand{\Q}{\mathbb Q}
\newcommand{\Z}{\mathbb Z}

\newcommand{\eps}{\varepsilon}

\newcommand{\Alt}{\mathrm{\bf A}}
\newcommand{\Sym}{\mathrm{\bf S}}
\newcommand{\Cyclic}{\mathrm{\bf C}}
\newcommand{\Dihedral}{\mathrm{\bf D}}
\newcommand{\M}{\mathrm{\bf M}}
\newcommand{\ASL}{\mathrm{\bf ASL}}
\newcommand{\AGL}{\mathrm{\bf AGL}}

\newcommand{\PSL}{\mathrm{\bf PSL}}
\newcommand{\PGL}{\mathrm{\bf PGL}}

\newcommand{\PYL}{\mathrm{\bf P\Gamma L}}

\usepackage{tikz}
\usetikzlibrary{matrix}

\newenvironment{xflalign*}
 {\setlength{\abovedisplayskip}{4pt}\setlength{\belowdisplayskip}{6pt}%
  \csname flalign*\endcsname}
 {\csname endflalign*\endcsname\ignorespacesafterend}
  
\usepackage{hyperref}
\hypersetup{
    colorlinks=true,
    linkcolor=blue,
    citecolor=blue,
    urlcolor=blue,
}

\newcommand{\subsectionsep}[1]{\bigskip\pagebreak[3]\hrule\bigskip\noindent{\large #1}\bigskip}
   
\title{
Belyi functions of the special weighted trees of $(2,3)$-type
}
\author{
Nikolai M. Adrianov\footnote{
Lomonosov Moscow State University, 119991 Moscow, Russia, e-mail: {\tt nadrianov@gmail.com}},
Alexander M. Vatuzov\footnote{
Lomonosov Moscow State University, 119991 Moscow, Russia, e-mail: {\tt amvatuzov@yandex.ru}}
}

\date{}

\begin{document}
\maketitle


%\begin{flushleft}
%УДК 519.17+511.238
%\end{flushleft}

\begin{abstract}
In this paper we calculate the Belyi functions of the weighted trees of $(2,3)$-type with primitive special edge rotation groups. There are 21 Galois orbits of these trees: 6 rational orbits, 12 orbits over quadratic fields, and three orbits over fields of degree 3, 4 and 6. The highest degree of the calculated Belyi function is 32. The calculations are performed using modular functions techniques.
\end{abstract}

\section*{Introduction}

A {\it dessin d'enfant} is a bicolored (connected) graph $\Gamma$ embedded into a compact oriented surface ${\cal X}$ in such a way that the complement ${\cal X}\setminus \Gamma$ is homeomorphic to a disjoint union of open discs. Let $X$ be a non-singular complex algebraic curve. A function $\beta$ on $X$ is called the {\it Belyi function} if $\beta$ has no critical values  other than $0$, $1$ and $\infty$. Then the preimage $\beta^{-1}[0, 1]$ is a dessin d'enfant on $X$. The inverse is also true: for an arbitrary dessin d'enfant $D=({\cal X}, \Gamma)$ there exists a corresponding Belyi pair $(X,\beta)$. Moreover the Belyi pair can be defined over a number field (an algebraic extension of $\Q$ of finite degree). A detailed introduction into the theory of dessins d'enfants can be found, e.g., in~\cite{LanZvo-2004}.

The action of the Galois group $Gal({\bar \Q}/\Q)$ on the Belyi pairs induces its action on the dessins d'enfants; the orbits of this action are finite. Therefore for any dessin d'enfant the stabilizer under this action is a finite index subgroup of $Gal({\bar \Q}/\Q)$; the correponding number field  is called {\it field of moduli} of the dessin. For a dessin of genus 0 we put $X=\PP^1$ and then its Belyi function $\beta\in{\bar \Q}(z)$ is uniquely determined up to a linear fractional substitution $z\to (az+b)/(cz+d)$. For all dessins in this paper the Belyi functions can be defined over field of moduli of the dessin; in a general case it is not always possible, see \cite{Cou-1994}.

The degree of a Belyi function is equal to the number of edges of the corresponding dessin. A triple of partitions $\pi(D)=(\lambda_0 \,|\, \lambda_1 \,|\, \lambda_2)$, where $\lambda_{0,1,2}\vdash n$ are the partitions constituted by the degrees of black vertices, white vertices and faces of a dessin $D$ with $n$ edges, is called the {\it passport} of the dessin. In terms of the Belyi pairs, passport describes ramifications of the Belyi function over $0$, $1$ and $\infty$. The passport of the dessin is an invariant of the Galois group action.

One can define a dessin $D$ by a pair of permutations of edges $(a,b)$ which correspond to counter-clockwise detours around black and white vertices. The cyclic structures of the permutations $a$, $b$ and $ab$ are the components of the passport of the dessin. A group generated by permutations $a$ and $b$ is called the {\it edge rotation group} and is denoted by $ER(D)$. The edge rotation group of a dessin is isomorphic to the monodromy group of the cover $\beta:X\to \PP^1(\CC)$; it is yet another invariant of the Galois group action.

\begin{figure}[!ht]
\centering
\includegraphics[angle=90,scale=0.8]{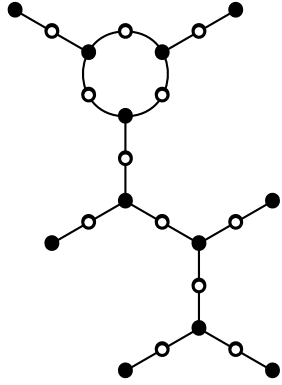}
\caption{``Alien'': a dessin with passport $(3^6 1^6 \,|\, 2^{12} \,|\, 21^1 3^1)$ and the edge rotation group $\M_{24}$.}
\label{fig:M24_alien}
\end{figure}

{\bf Example.} The ``alien'' dessin shown in fig.~\ref{fig:M24_alien} is a dessin d'enfant with 24 edges and the passport $(3^6 1^6 \,|\, 2^{12} \,|\, 21^1 3^1)$; its edge rotation group is the marvelous Mathieu group $\M_{24}$, see \cite{Zvo-1998}. The Belyi function of this dessin defined over $\Q(\sqrt{-7})$ is calculated in \cite{Vatuzov-2021}.
 
A dessin of genus $0$ with only one face of degree distinct from 1 is called the {\it weighted tree}.  Weighted trees can be understood combinatorially as plane trees in which some edges are replaced by multiple ones, see~\cite{AdrPakZvo-2020}, \cite{PakZvo-2014}. In \cite{AdrZvo-2013} weighted trees with special primitive edge rotation groups\footnote{A primitive permutation group on $n$ points is called {\it special} if it is different from $\Cyclic_n$,  $\Dihedral_n$, $\Alt_n$ and $\Sym_n$.} are classified. There are 184 such trees, which are subdivided into (at least) 85 Galois orbits and generate 34 primitive groups (the highest degree is 32).

We say that a dessin d'enfant is of {\it $(2,3)$-type} if all its black vertices are of degree 1 or 3 and all white vertices are of degree 1 or 2. In other words, the permutations $a$ and $b$ are of degree $3$ and $2$ respectively. In this paper we calculate the Belyi functions for all weighted trees of  $(2,3)$-type with special primitive edge rotation groups. There are 43 such trees in 21 Galois orbits; the number of edges varies from 6 up to 32. The Belyi functions for some of them were known before; we used them to test our calculation method. The answers are included here for completeness.

As far as we know, the Belyi functions for orbits 12.4, 12.13, 13.1, 14.2, 17.1, 20.1, 24.2, 32.1 are published for the first time (in \cite{AdrPakZvo-2020} the fields of definition for these orbits are given, but not the Belyi functions itself).

\section*{Method of calculation}

The Belyi function $\beta$ of a dessin of genus 0 is a rational function in one variable. One can use the passport write down the numerators and denominators of $\beta$ and $\beta-1$ as products of polynomials with indeterminate variables. For example, the passport of a weighted tree  of $(2,3)$-type with $n$ edges has the form $(3^{p_3} 1^{p_1} \,|\, 2^{q_2} 1^{q_1} \,|\, (n-r)^1\, 1^r)$; put the pole of degree $n-r$ at $z=\infty$ and then the Belyi function takes the form
\begin{equation}
\label{eq:system}
\beta = \frac{P_3^3\cdot P_1}{cR} = \frac{Q_2^2\cdot Q_1}{cR} + 1,
\end{equation}
where $c\in \CC$ is a constant, $P_3$, $P_1$, $Q_2$, $Q_3$, $R\in \CC[z]$ are monic polynomials without multiple roots, the degrees of the polynomials are the multiplicities of the passport: $\deg P_i = p_i$, $\deg Q_i = q_i$, $\deg R = r$.

The equation (\ref{eq:system}) gives us a system of non-linear algebraic equations in the coeficients of the polynomials $P_i$, $Q_i$, $R$. The system is underdetermined: there is still a freedom of affine substitution $z\to Az+B$. For $(2,3)$-dessins that appear in the current paper we can add the following two normalization conditions:
\begin{equation}
\label{eq:system2}
\begin{array}{l}
\text{the sum of coordinates of the black vertices of degree 3} = 0,\\
\text{the sum of coordinates of the white vertices of degree 2} = 1.
\end{array}
\end{equation}
These conditions set the values of two variables in the system (\ref{eq:system}):
\begin{equation}
\label{eq:system3}
\mathrm{coeff}_{p_3-1} P_3 = 0,
\quad
\mathrm{coeff}_{q_2-1} Q_2 = -1.
\end{equation}

The most straightforward approach here is to solve the resulting system using the Gr\"obner bases approach. Unfortunately the complexity of this methods grows rapidly as the number of edges increases; it is practically not applicable for dessins with more than 12 edges.

A more promising approach uses approximate numerical calculations. The general calculation scheme is as follows. We solve the system obtained from (\ref{eq:system}) and (\ref{eq:system2}) numerically with high enough precision and then express coefficients as algebraic numbers using LLL-algorithm\footnote{
The abbreviation LLL is made up of the first letters of the names of the authors A. K. Lenstra, H. W. Lenstra, and L. Lov\'asz. In SageMath the algorithm is implemented by {\tt algdep} function.}.
This scheme was first used to calculate a Belyi function in \cite{Cou-1994}. The problem here is to find a good enough starting point such that the Newton method converges to the Belyi function of our dessin\footnote{Recall that we write down the system using only the passport of the dessin. Meanwhile the number of the dessins with the given passport can be huge even for a relatively small number of edges.}. We use the modular functions method based on ideas from \cite{Mon-2014} and \cite{SijVoi-2013}.

Let $D$ be a $(2,3)$-dessin with $n$ edges defined by permutations $a$, $b\in \Sym_n$. Then the edge rotation group $ER(D)$ can be presented as a factor of the modular group
$$
\Gamma = \PSL_2(\Z) \simeq \langle S,T \mid S^2=(ST)^3=1 \rangle
$$
by a subgroup $\Gamma_D\subset \Gamma$ of index $n$. The group $\Gamma$ acts on the completed upper half-plane $\HH^* = \HH\cup\Q\cup\{i\infty\}$ defining the well-known tesselation of $\HH^*$ by ideal triangles. The fundamental domain of $\Gamma_D$ consists of $n$ such triangles; it can be easily calculated from $a$ and $b$. We get the commutative triangle
%%
%% TODO: j(i) = 1728, j(e^{2\pi i/3}) = 0
%%
\begin{center}
\vskip -0.8em
\begin{tikzpicture}
\matrix (m) [matrix of math nodes, row sep=1em, column sep=0.6em]
{
& \HH &\\
& & \\
X & & \,\,\PP^1 \\
};
  \path[-stealth]
    (m-1-2) edge node [left] {$\scriptstyle /\Gamma_D$} (m-3-1) % /\Gamma_D = t_D
			edge node [right] {$\scriptstyle \,J$} (m-3-3) % /\Gamma = j
    (m-3-1) edge node [above] {$\scriptstyle \;\, \beta$} (m-3-3);
\end{tikzpicture}
\end{center}
\vskip -1.2em
Here $J(\tau) = \frac{1}{1728} j(\tau)$, where $j(\tau)$ is a classical modular function invariant under the action of $\Gamma$. For a dessin $D$ of genus $0$, we have $X\simeq\PP^1$ and the factorization by $\Gamma_D$ is defined by a modular function $t_D:\,\HH\to \PP^1$ invariant under the action of $\Gamma_D$. 

Let $m$ be the degree of the pole of $\beta$ at $z=\infty = t_D(i\infty)$. Then the function $t_D$ is invariant under substitution $\tau\to \tau+m$ and can be expressed via $\zeta = q^{1/m} = e^{2\pi i\tau/m}$ and $z=t_D$ has a pole of order 1 at $\zeta=0$. Then the power series expansion of $t_D$ in variable $\zeta$ takes the form
$$
t_D(\zeta) = \frac{c_{-1}}{\zeta} + c_0 + \sum_{k=1}^\infty c_k \zeta^k, \quad c_{-1}\ne 0.
$$
Denote by $[t_D]_N(\zeta)$ a partial sum of this expansion:
$$
[t_D]_N(\zeta) = \frac{c_{-1}}{\zeta} + c_0 + \sum_{k=1}^N c_k \zeta^k.
$$
Choose $N$ points $\{x_i\}_{i=1}^N$ at the boundary of the fundamental domain of $\Gamma_D$, putting several points for each arc bounding the fundamental domain. Write down the conditions $[t_D]_{N}(x_i)=[t_D]_{N}(M_i(x_i))]$, where $M_i\in\Gamma_D$ is a transformation that defines the pasting of the arc containing $x_i$. Add the conditions ({\ref{eq:system2}) written down in terms of $[t_D]_N$. We get a system of linear equations of size $(N+2)\times (N+2)$; solving it we obtain approximate values for the coefficients $\{c_{-1},c_0,\ldots c_{N}\}$. Substituting the coordinates of the vertices of the triangles of the fundamental domain of $\Gamma_D$ into the approximation of $[t_D]_N$, we obtain the approximate coordinates for the vertices and of the dessin d'enfant, and thus, the approximate values for coefficients in (\ref{eq:system}). We use these values as a starting point for the Newton method.

Applying the method described above we have calculated (in almost automatic mode) the Belyi functions of weighted trees of degree up to 32. The method is applicable also to the dessins d'enfants which are not of $(2,3)$-type: any dessin can be reduced to a triangulation and then to a dessin of $(2,3)$-type. However the procedure can increase the number of edges by 6 times, and hence, increase the complexity of the system of algebraic equations. In future works we are going to calculate more Belyi functions and to determine the limits of the method.

At the time of writing, the highest records are the calculations of the Belyi functions with monodromy groups $J_2$ and $Co_3$ of degree 100 and 276 respectively, see \cite{Mon-2017}, \cite{Mon-2018}. Similar methods were used in \cite{MSSV-2018} to compute the database of Belyi maps (of degree up to 9).

\section*{Results of calculation}

In this section we present the list of calculated Belyi functions. We try to write down the Belyi functions in a most compact way (orienting to the amazing classical formulas by F.\,Klein). 

%Writing down the expression for the Belyi function (\ref{eq:system}), 
In (\ref{eq:system})
we placed the only pole of degree $>1$ to the point $z=\infty$. This leaves us a freedom of affine substitutions $z\to Az+B$. For a weighted tree with only one face of degree 1 we place the second pole at the point $z=0$. In all other cases we take the second leading coefficient of the polynomial $P_3$, $P_1$, $Q_2$, $Q_1$ or of the numerator of the Belyi function and assign it to zero. Our choice of the polynomial here is more or less arbitrary.

There is still a freedom of multiplicative substitution $z\to Az$. We choose the value of $A$ in such a way that
\begin{enumerate}
\item[(1)] the leading coefficients of $P_3$, $P_1$, $Q_2$, $Q_1$ and $R$ are equal to $1$ (in several cases we get more compact formulas if we allow the polynomial $R$ to be non-monic);
\item[(2)] all other coefficients are integer algebraic numbers (the smallest possible with respect to the norm).
\end{enumerate}
We get into a trouble with this procedure only for the orbit 32.1 when the coefficients have non-unique factorization. Is there a proper way to construct a minimal model of the Belyi function?  It is an interesting question.

%%%%%%%%%%%%%%%%%%%%%%%%%%%%%%%%%%%%%%%%%%

\subsectionsep{$n=6$, $(3^2 \,|\, 2^2 1^2 \,|\, 5^1 1^1)$, $ER=\PSL_2(5)$}

A weighted tree (orbit 6.1, see.~\cite{AdrZvo-2013}), the orbit consists of a single dessin, see fig.~\ref{fig:PSL_2_5}.
The Belyi function is the covering map $X_0(5)\to X(1)$ of the classical modular curves; it was actually found in a classical work of  H.\,Schwarz\cite{Schwarz-1873} 150 years ago (the famous Schwarz's list of algebraic hypergeometric functions). Few years later F.\,Klein depicted dessins d'enfants of $X_0(N)$ for $N=2,3,4,5,7,13$ and calculated their Belyi functions, see~\cite{Klein-1878}, p.~143 and figures 10, 12, 14 at pp.~136--138.

\begin{figure}[!ht]
\centering
\begin{subfigure}[b]{0.45\textwidth}
\centering
\includegraphics[scale=0.8]{images/wt-601.mps}
\vspace{2cm}
\end{subfigure}
\begin{subfigure}[b]{0.45\textwidth}
\centering
\includegraphics[width=5cm]{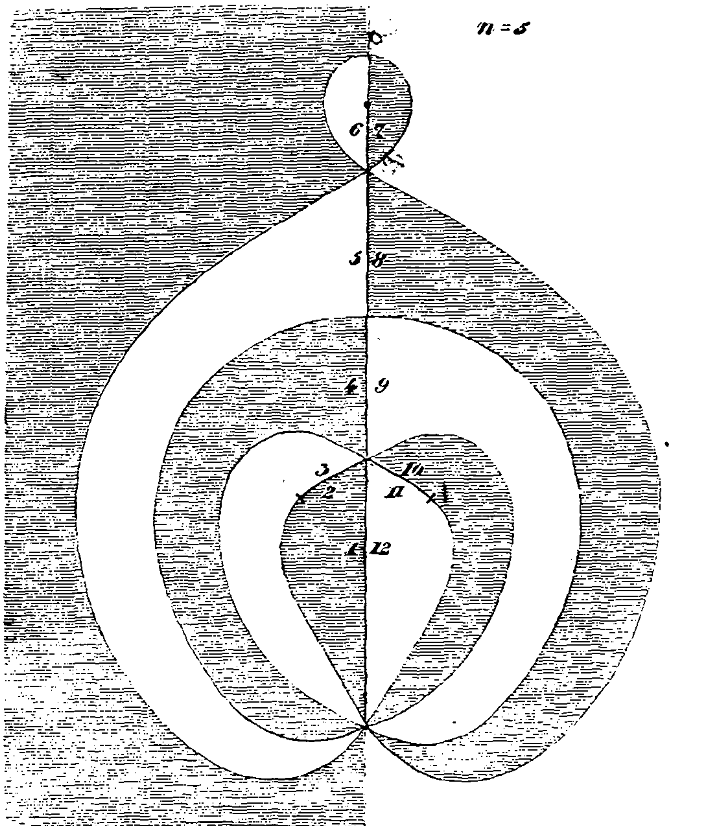}
\end{subfigure}
\caption{Group $\PSL_2(5)$: orbit 6.1 and a drawing by F.\,Klein.}
\label{fig:PSL_2_5}
\end{figure}
%$$
%\beta = \frac{(z^2-5)^3}{27(2z-5)}
%$$
%$$
%\beta -1 = \frac{(z^2-3z+1)^2(z^2+6z+10)}{27(2z-5)}
%$$
$$
\beta = \frac{(z^2+10z+5)^3}{1728z} =
\frac{(z^2+4z-1)^2(z^2+22z+125)}{1728z}
+ 1.
$$

%%%%%%%%%%%%%%%%%%%%%%%%%%%%%%%%%%%%%%%%%%

\subsectionsep{$n=7$, $(3^2 1^1 \,|\, 2^3 1^1 \,|\, 6^1 1^1)$, $ER=\AGL_1(7)$}

A weighted tree (orbit 7.1, see~\cite{AdrZvo-2013}), the orbit consists of two dessin, the one shown in fig.~\ref{fig:AGL_1_7} and his mirror image. The dessin is also used in \cite{SijVoi-2013} to demonstrate the method of calculation based on the modular functions techniques.

\begin{figure}[!ht]
\centering
\par\vspace{0.3cm}
\includegraphics[scale=0.8]{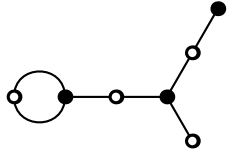}
\caption{Group $\AGL_1(7)$: orbit 7.1 of size 2.}
\label{fig:AGL_1_7}
\end{figure}

$$
\beta = \frac{P_3^3 P_1}{cz} = \frac{Q_2^2\cdot Q_1}{cz} + 1,
$$
where
\begin{xflalign*}
c = -16 (1 - \omega) (1 + \omega)^3 (3 - \omega)^7
\qquad\qquad\qquad\qquad
\omega = \frac{1\pm\sqrt{-3}}{2}
\end{xflalign*}
\begin{xflalign*}
P_3 ={}
  & z^2 - 2 \omega (1 + \omega) (4 - \omega) z - 4 (2 - \omega)
  \qquad\qquad\quad\,
  P_1 = z - 18\omega
\end{xflalign*}
\begin{xflalign*}
Q_2 = {}
  & z^3 - 2\omega (1 + \omega) (7 - 3\omega) z^2 + 4\omega (1 + \omega)^3 (4 - \omega) z - 12 (1 + \omega)
\end{xflalign*}
\begin{xflalign*}
Q_1 = {}
  & z + 8 (1 - 2\omega).
\end{xflalign*}

%%%%%%%%%%%%%%%%%%%%%%%%%%%%%%%%%%%%%%%%%%

\subsectionsep{$n=7$, $(3^2 1^1 \,|\, 2^2 1^3 \,|\, 7^1)$, $ER=\PSL_3(2)$}

An ordinary tree (orbit 7.2, see~\cite{AdrZvo-2013}), the orbit consists of two dessins, see fig.~\ref{fig:PSL_3_2}. The corresponding Belyi function is a polynomial (the so-called {\it Shabat polynomial}) and calculated by F.\,Klein \cite{Klein-1878b}, p.426-427, and later in \cite{BePeZvo-1992}.

\begin{figure}[!ht]
\centering
\begin{subfigure}[b]{0.45\textwidth}
\centering
\includegraphics[scale=0.8]{images/wt-702.mps}
\vspace{1.5cm}
\end{subfigure}
\begin{subfigure}[b]{0.45\textwidth}
\centering
\includegraphics[width=5cm]{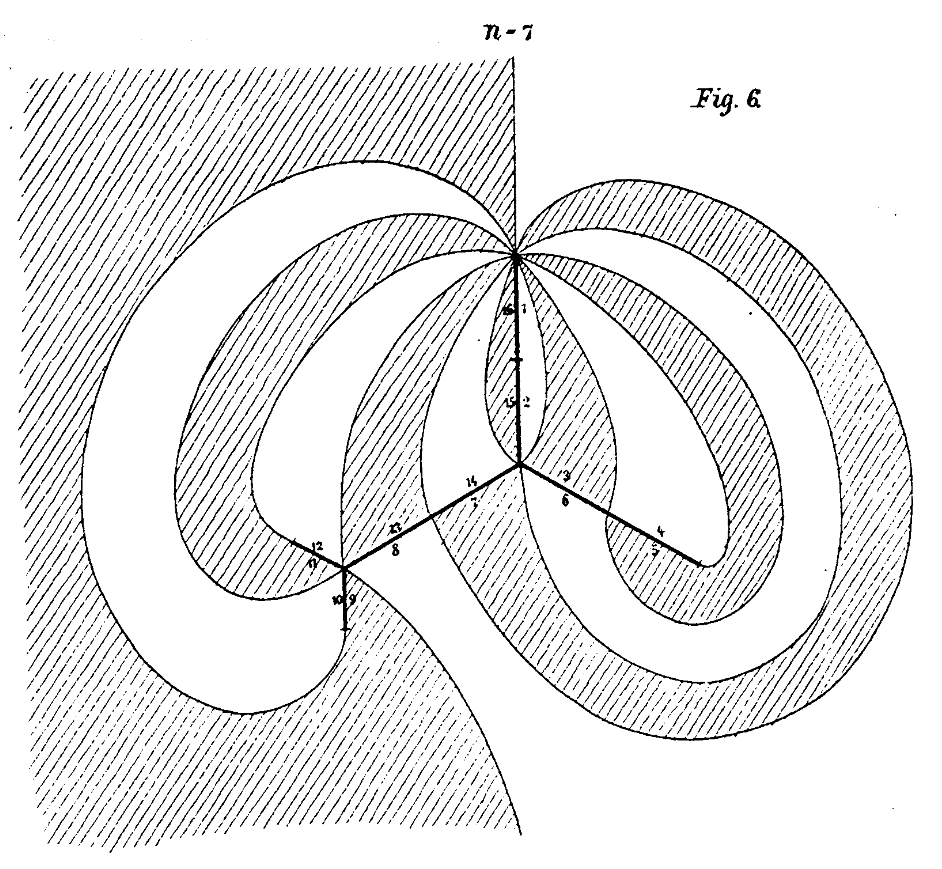}
\end{subfigure}
\caption{Group $\PSL_3(2)$: orbit 7.2 of size 2 and a drawing by F.\,Klein.}
\label{fig:PSL_3_2}
\end{figure}

$$
\beta = \frac{P_3^3\cdot P_1}{c} = \frac{Q_2^2\cdot Q_1}{c} + 1,
$$
where
\begin{xflalign*}
c = {} & 2^6 \cdot 3^3\cdot \omega^{14}
\qquad\qquad\qquad
\qquad\qquad\qquad
\omega = \frac{1\pm\sqrt{-7}}{2}\\
P_3 = {} & z^2 + 28
\qquad\quad
P_1 = z + 7\omega
\qquad\quad
Q_2 = z^2 + 6\omega z + 4\\[4pt]
Q_1 ={} &
  \bigl(z^2 - 4z - 4(1-\omega)(1-6\omega)\bigr)
  \bigl(z - \omega (3+2\omega)\bigr).
\end{xflalign*}
%\beta_m1 = Q2^2* Q1 / c

%%%%%%%%%%%%%%%%%%%%%%%%%%%%%%%%%%%%%%%%%%

\subsectionsep{$n=8$, $(3^2 1^2 \,|\, 2^4 \,|\, 7^1 1^1)$, $ER=\PSL_2(7)$}

A weighted tree (orbit 8.8, see~\cite{AdrZvo-2013}), the orbit consists of one dessin, see fig.~\ref{fig:PSL_2_7}.
The Belyi function is the covering map $X_0(7)\to X(1)$ of the classical modular curves; it was calculated by F.\,Klein~\cite{Klein-1878}.
\begin{figure}[!ht]
\centering
\begin{subfigure}[b]{0.45\textwidth}
\centering
\includegraphics[scale=0.8]{images/wt-809.mps}
\vspace{1.5cm}
\end{subfigure}
\begin{subfigure}[b]{0.45\textwidth}
\centering
\includegraphics[width=5.5cm]{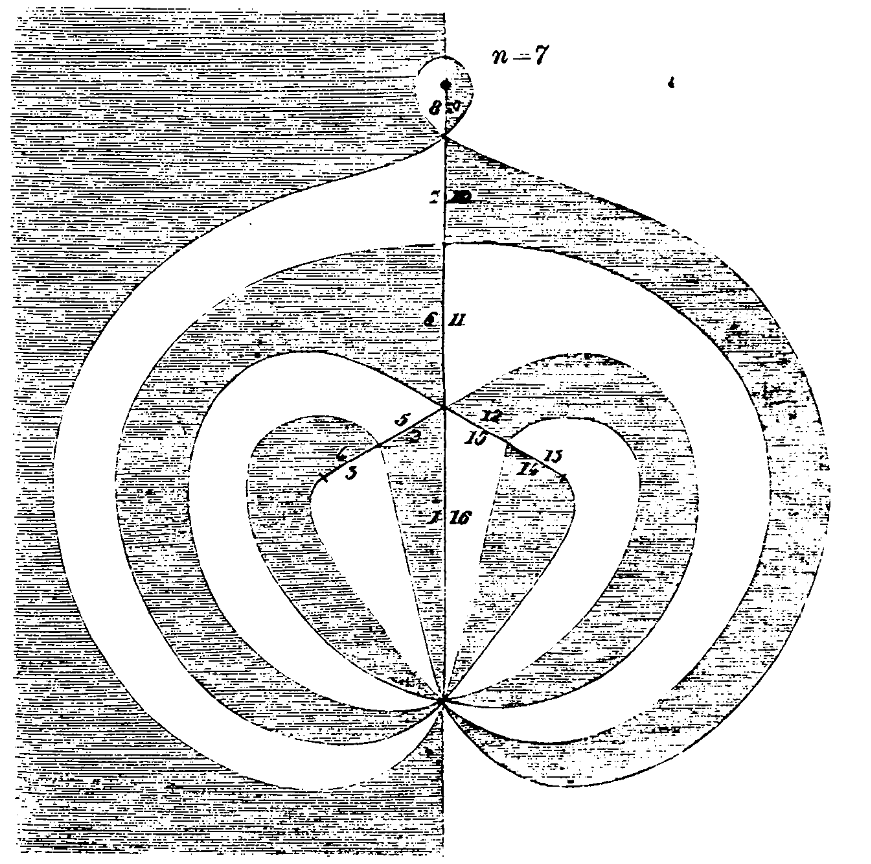}
\end{subfigure}
\caption{Group $\PSL_2(7)$: orbit 8.8 of size 1 and a drawing by F.\,Klein.}
\label{fig:PSL_2_7}
\end{figure}

$$
\beta = \frac{(z^2 + 5z + 1)^3(z^2 + 13z + 49)}{1728z} = 
\frac{(z^4 + 14z^3 + 63z^2 + 70z - 7)^2}{1728z} + 1.
$$

%%%%%%%%%%%%%%%%%%%%%%%%%%%%%%%%%%%%%%%%%%

\subsectionsep{$n=8$, $(3^2 1^2 \,|\, 2^3 1^2 \,|\, 8^1)$, $ER=\PGL_2(7)$}

An ordinary tree (orbit 8.15, see~\cite{AdrZvo-2013}), the orbit consists of one dessin, see fig.~\ref{fig:PGL_2_7}.
The corresponding Belyi function (Shabat polynomial) is calculated in~\cite{BePeZvo-1992}.

\begin{figure}[!ht]
\centering
\begin{subfigure}[b]{\textwidth}
\centering
\begin{subfigure}[b]{0.35\textwidth}
\centering
\includegraphics[angle=90,scale=0.8]{images/wt-816.mps}
\end{subfigure}
\begin{subfigure}[b]{0.35\textwidth}
\centering
\includegraphics[scale=0.8]{images/wt-817.mps}
\vspace{0.21cm}
\end{subfigure}
\end{subfigure}
%\caption{Group $\PGL_2(7)$: orbit 8.15 of size 2.}
\caption{Группа $\PGL_2(7)$: орбита 8.15 размера 2.}
\label{fig:PGL_2_7}
\end{figure}

$$
c = 2^4 \cdot 3^3 \cdot \omega^7 (2+\omega)
\qquad\qquad
\omega = 1\pm 2\sqrt{2}
$$
$$
\beta = \frac{1}{c}
\left(z^2 + \omega^3\right)^3 \left(z^2 - 8z + \omega^2 (8+\omega)\right)
$$
$$
\beta - 1 = \frac{1}{c}
\left(z^3 - 7z^2 + \omega^2 (6+\omega) z - \omega^3\right)^2 \left(z^2 + 6z + 7(5+2\omega)\right).
$$

%%%%%%%%%%%%%%%%%%%%%%%%%%%%%%%%%%%%%%%%%%

\subsectionsep{$n=9$, $(3^3 \,|\, 2^3 1^3  \,|\, 8^1 1^1)$, $ER=\AGL_2(3)$}

A weighted tree (orbit 9.2, see~\cite{AdrZvo-2013}), the orbit consists of two dessins, see fig.~\ref{fig:AGL_2_3}.
The Belyi function of the dessin can be found in the database of Belyi pairs~\cite{MSSV-2018}.

\begin{figure}[!ht]
\centering
\includegraphics[scale=0.8]{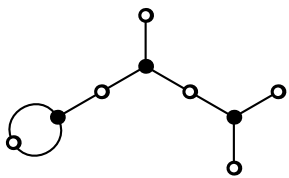}
\caption{Group $\AGL_2(3)$: orbit 9.2 of size 2.}
\label{fig:AGL_2_3}
\end{figure}

$$
\beta = \frac{P_3^3}{cz} = \frac{Q_2^2\cdot Q_1}{cz} + 1,
$$
where
$$
c = -3^3 \cdot \omega^{12}
\qquad\qquad\qquad
\omega = 1 \pm \sqrt{-2}
$$
\begin{xflalign*}
P_3 = {}
  & z^3 + 2^3 \omega z^2 - 2^2 \omega (1-4\omega) z - 2^3 \\
Q_2 = {}
  & z^3 + 5 \omega z^2 - \omega (1-4\omega) z + 1 \\
Q_1 = {}
  & z^3 + 2\cdot 7\omega z^2 - \omega (10-67\omega) z - 2^9.
\end{xflalign*}

%%%%%%%%%%%%%%%%%%%%%%%%%%%%%%%%%%%%%%%%%%

\subsectionsep{$n=9$, $(3^3 \,|\, 2^4 1^1  \,|\, 7^1 1^2)$, $ER=\PSL_2(8)$}

A weighted tree (orbit 9.4, see~\cite{AdrZvo-2013}), the orbit consists of one dessin, see fig.~\ref{fig:PSL_2_8}.
The dessin is the only realization of its passport. The Belyi function is calculated in \cite{PakZvo-2018} (sporadic unitree $L$).

\begin{figure}[!ht]
\centering
\includegraphics[angle=-90,scale=0.8]{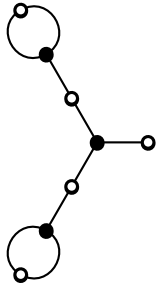}
\caption{Group $\PSL_2(8)$: orbit 9.4 of size 1.}
\label{fig:PSL_2_8}
\end{figure}

%$$
%\beta = -\frac{(z^3 + 7 z^2 + 315 z + 2037)^3}{2^{21} \cdot 3^3(z^2 + 343)}
%$$
%$$
%\beta - 1 = -\frac{(z^4 + 4 z^3 + 486 z^2 + 1476 z + 46305)^2 (z + 13)}{2^{21} \cdot 3^3(z^2 + 343)}.
%$$
$$
\beta = -\frac{(z^3 - 16z^2 + 160z - 384)^3}{2^{14}\cdot 3^3 (z^2 - 13z + 128)}
$$
$$
\beta - 1 = -\frac{z(z^4 - 24z^3 + 336z^2 - 2240z + 8064)^2}{2^{14}\cdot 3^3 (z^2 - 13z + 128)} + 1
$$
or (if we do not require the polynomial $R$ to be monic) the Belyi function can be written as
$$
\beta = -\frac{(z^3 - 8z^2 + 40z - 48)^3}{1728 (2z^2 -13z + 64)}
$$
$$
\beta - 1 = -\frac{z(z^4 - 12z^3 + 84z^2 - 280z + 504)^2}{1728 (2z^2 -13z + 64)}.
$$

%%%%%%%%%%%%%%%%%%%%%%%%%%%%%%%%%%%%%%%%%%

\subsectionsep{$n=9$, $(3^2 1^3 \,|\, 2^4 1^1 \,|\, 9^1)$, $ER=\PYL_2(8)$}

An ordinary tree (orbit 9.6, see~\cite{AdrZvo-2013}), the orbit consists of two dessins, see fig.~\ref{fig:PYL_2_8}.
The corresponding Shabat polynomial can be found in the catalog of plane trees with 9 edges \cite{Koch-2007}.

\begin{figure}[!ht]
\centering
\includegraphics[scale=0.8]{images/wt-907.mps}
\caption{Group $\PYL_2(8)$: orbit 9.6 of size 2.}
\label{fig:PYL_2_8}
\end{figure}

$$
\beta = \frac{P_3^3\cdot P_1}{c} = \frac{Q_2^2\cdot Q_1}{c},
$$
where
$$
c = -2^{15} (1-\omega) (3-\omega)^7
\qquad\qquad\qquad
\omega = \frac{1 \pm \sqrt{-3}}{2}
$$
\begin{xflalign*}
P_3 = {}
 & z^2 - \omega (1+\omega) (3-\omega)^3 \\[4pt]
P_1 = {}
 & z^3 - 9 z^2 - (1+\omega)^3 (17-6\omega) z + \omega (1+\omega)^2 (3-\omega)^2 (47-12\omega) \\[4pt]
Q_2 = {}
 & z^4 - 8 z^3 - 6\omega (3-\omega) (3+2\omega) z^2 + {} \\
 & {} + 2^3\cdot 7 \omega (3-\omega) (3+2\omega) z - (3-\omega)^3 (17-6\omega) \\[4pt]
Q_1 = {} & z + 7.
\end{xflalign*}

%%%%%%%%%%%%%%%%%%%%%%%%%%%%%%%%%%%%%%%%%%

\subsectionsep{$n=10$, $(3^3 1^1 \,|\, 2^5 \,|\, 8^1 1^2)$, $ER=\PGL_2(9)$}

A weighted tree (orbit 10.1, see~\cite{AdrZvo-2013}), the orbit consists of one dessin, see fig.~\ref{fig:PGL_2_9}.
The dessin is the only realization of its passport. The Belyi function is calculated in \cite{PakZvo-2018} (sporadic unitree $M$).

\begin{figure}[!ht]
\centering
\includegraphics[angle=-90,scale=0.8]{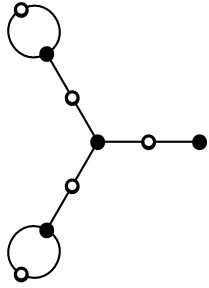}
\caption{Group $\PGL_2(9)$: orbit 10.1 of size 1.}
\label{fig:PGL_2_9}
\end{figure}

%$$
%\beta = -\frac{4(z^3 + 3z - 2)^3 (z - 2)}{81\left(z^2 + \frac{2}{3}z + \frac{11}{3}\right)}
%$$
%$$
%\beta - 1 = -\frac{4\left(z^5 - z^4 + 4z^3 - 8z^2 + \frac{7}{2}z - \frac{19}{2}\right)^2}{81\left(z^2 + \frac{2}{3}z + \frac{11}{3}\right)}
%$$
$$
\beta = -\frac{z(z^3 + 36z^2 + 540z + 2592)^3}{2^6\cdot 3^{12} (z^2 + 28z + 324)}
$$
$$
\beta - 1 = -\frac{(z^5 + 54z^4 + 1296z^3 + 15552z^2 + 87480z + 104976)^2}{2^6\cdot 3^{12} (z^2 + 28z + 324)}
$$
or (if we do not require the polynomial $R$ to be monic) the Belyi function can be written as
$$
\beta = -\frac{z(z^3 + 12z^2 + 60z + 96)^3}{1728 (3z^2 + 28z + 108)}
$$
$$
\beta - 1 = -\frac{(z^5 + 18z^4 + 144z^3 + 576z^2 + 1080z + 432)^2}{1728 (3z^2 + 28z + 108)} + 1.
$$

%%%%%%%%%%%%%%%%%%%%%%%%%%%%%%%%%%%%%%%%%%

\subsectionsep{$n=11$, $(3^3 1^2 \,|\, 2^4 1^3 \,|\, 11^1)$, $ER=\PSL_2(11)$}

An ordinary tree (orbit 11.1, see~\cite{AdrZvo-2013}), the orbit consists of two dessins, see fig.~\ref{fig:PSL_2_11}.
The Shabat polynomial is calculated by F.\,Klein, see~\cite{Klein-1879}, p.~547 (the drawing is taken from \cite{Klein-1878b}), and later in \cite{CasCou-1999} and \cite{JonZvo-2021}.

\begin{figure}[!ht]
\centering
\begin{subfigure}[b]{0.45\textwidth}
\centering
\includegraphics[scale=0.8]{images/wt-1101.mps}
\vspace{1cm}
\end{subfigure}
\begin{subfigure}[b]{0.45\textwidth}
\centering
\includegraphics[width=6cm]{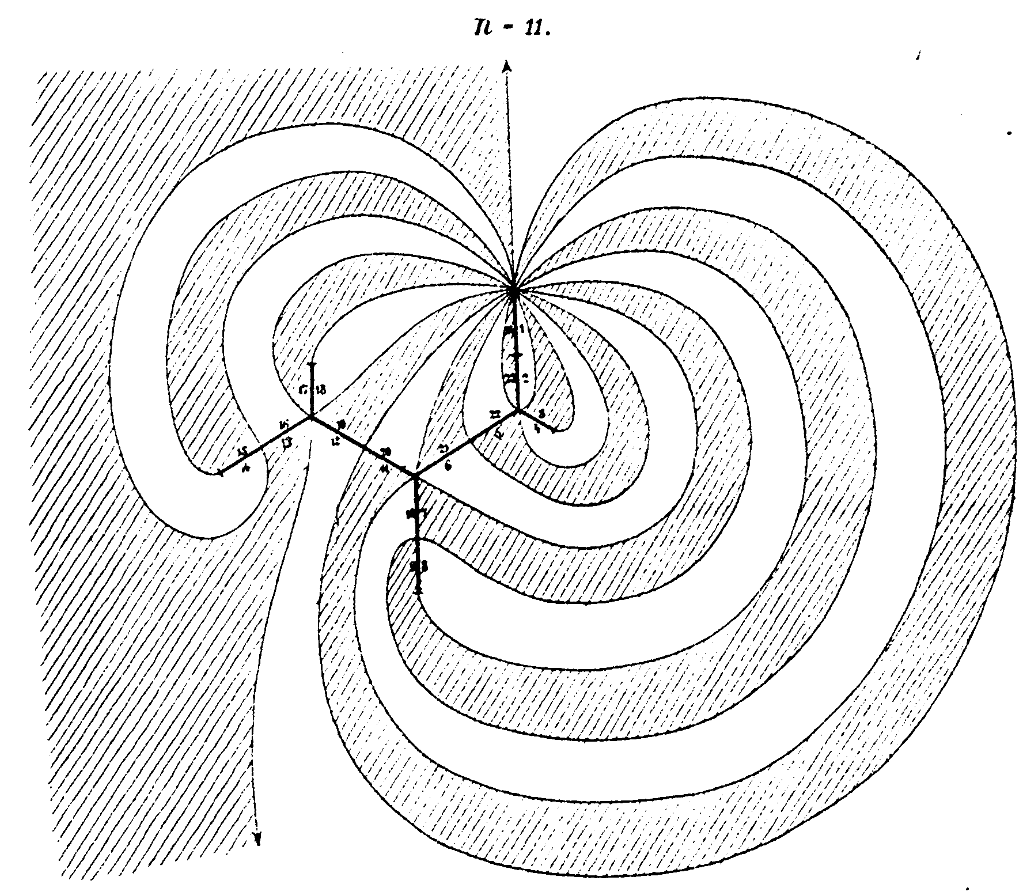}
\end{subfigure}
\caption{Group $\PSL_2(11)$ acting on 11 points: orbit 11.1 of size 2 and a drawing by F.\,Klein.}
\label{fig:PSL_2_11}
\end{figure}

$$
\beta = \frac{P_3^3\cdot P_1}{c} = \frac{Q_2^2\cdot Q_1}{c},
$$
where
\begin{xflalign*}
c = {} & -1728
  \qquad\qquad\qquad\qquad
  \qquad\qquad\qquad\qquad
  \omega=\frac{1\pm\sqrt{-11}}{2} \\
P_3 = {}
  & \bigl(z - (1+\omega)\bigr) \bigl(z^2 + (2+\omega) z - (1-\omega)\bigr) \\
P_1 = {}
  & z^2 - 3z - 2\omega^2 \\
Q_2 = {}
  & \bigl(z - (1-\omega)\bigr) \bigl(z^3 - (1+\omega) z^2 - (1+2\omega) z - 3 (1-\omega)\bigr) \\
Q_1 = {}
  & z^3 + 4z^2 + \omega^4 z + 2 (5-6\omega).
\end{xflalign*}

%%%%%%%%%%%%%%%%%%%%%%%%%%%%%%%%%%%%%%%%%%

\subsectionsep{$n=12$, $(3^4 \,|\, 2^4 1^4 \,|\, 11^1 1^1)$, $ER=\M_{12}$}

A weighted tree (orbit 12.4, see~\cite{AdrZvo-2013}), the orbit consists of two dessins, see fig.~\ref{fig:M_12a}. 

\begin{figure}[!ht]
\centering
\includegraphics[scale=0.8]{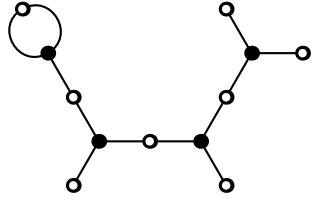}
\caption{Group $\M_{12}$: orbit 12.4 of size 2.}
\label{fig:M_12a}
\end{figure}

$$
\beta = \frac{P_3^3}{cz} = \frac{Q_2^2\cdot Q_1}{cz},
$$
where
\begin{xflalign*}
c = {} & 2^6 \cdot 3^{15}
  \qquad\qquad\qquad\qquad
  \qquad\qquad\qquad\qquad
  \qquad\qquad
  \omega=\frac{1\pm\sqrt{-11}}{2} \\
P_3 = {}
  & z^4 + 2\cdot 11\omega z^3 - 3\cdot 11 (13-5\omega) z^2 - 2\cdot 11 (56+15\omega) z - 11\omega^6 \\
Q_2 = {}
  & z^4 + 2^4\omega z^3 -  3\cdot 5 (13-5\omega) z^2 - 2^2 (56+15\omega) z + \omega^6 \\
Q_1 = {}
  & z^4 + 2\cdot 17\omega z^3 - 3 (2-\omega) (172+21\omega) z^2 - 2 (4216+1515\omega) z - 11^3\omega^6.
\end{xflalign*}

%%%%%%%%%%%%%%%%%%%%%%%%%%%%%%%%%%%%%%%%%%

\subsectionsep{$n=12$, $(3^3 1^3 \,|\, 2^6 \,|\, 11^1 1^1)$, $ER=\M_{12}$}

A weighted tree (orbit 12.5, see~\cite{AdrZvo-2013}), the orbit consists of two dessins, see fig.~\ref{fig:M_12b}.
The dessin was first presented by G.\,Malle at the Luminy conference in 1993 and was popularized by A.\,Zvonkin.
The dessin was given the name ``monsieur Mathieu''\footnote{L. Le Bruyn, \url{http://www.neverendingbooks.org/monsieur-mathieu}}.
The Belyi function is calculated by N.\,Magot, see~\cite{Zvo-1998}.

\begin{figure}[!ht]
\centering
\includegraphics[scale=0.8]{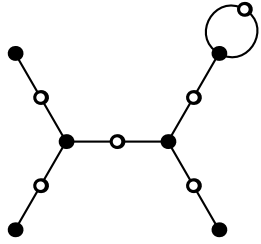}
\caption{Group $\M_{12}$: orbit 12.5 of size 2.}
\label{fig:M_12b}
\end{figure}

$$
\beta = \frac{P_3^3\cdot P_1}{cz} = \frac{Q_2^2}{cz},
$$
where
\begin{xflalign*}
c = {} & 2^6 \cdot 3^{12}
  \qquad\qquad\qquad\qquad
  \qquad\qquad\qquad\qquad
  \qquad\qquad
  \omega=\frac{1\pm\sqrt{-11}}{2} \\
P_3 = {}
  & z^3 - 3^2 (1-\omega) z^2 + 3 \omega^2 (2+3\omega) z + \omega^4 \\
P_1 = {}
  & z^3 - 17 (1-\omega) z^2 - (1-\omega) (1+\omega) (29+35\omega) z + 11^2 \omega^4 \\
Q_2 = {}
  & z^6 - 2\cdot 11(1-\omega) z^5 + 3\cdot 11 \omega^2(1+3\omega) z^4 +
     2\cdot 11\omega^4 (17+3\omega) z^3 - {} \\
  & \qquad\qquad\qquad
     {} - 3\cdot 11 \omega^4 (25-23\omega) z^2 +
     2\cdot 3\cdot 11\omega^6 (2+3\omega) z - 11\omega^8.
\end{xflalign*}

%%%%%%%%%%%%%%%%%%%%%%%%%%%%%%%%%%%%%%%%%%

\subsectionsep{$n=12$, $(3^4 \,|\, 2^5 1^2 \,|\, 10^1 1^2)$, $ER=\PGL_2(11)$}

A weighted tree (orbit 12.13, see~\cite{AdrZvo-2013}), the orbit consists of two dessins, see fig.~\ref{fig:PGL_2_11}.

\begin{figure}[!ht]
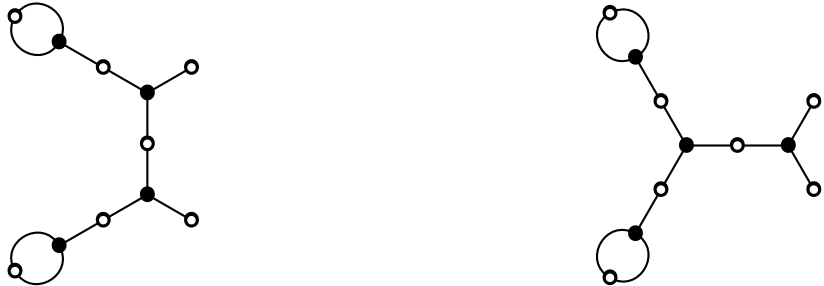

\centering
\begin{subfigure}[b]{\textwidth}
\centering
\begin{subfigure}[b]{0.45\textwidth}
\centering
\includegraphics[angle=90,scale=0.8]{images/wt-1216.mps}
\end{subfigure}
\begin{subfigure}[b]{0.45\textwidth}
\centering
\includegraphics[angle=90,scale=0.8]{images/wt-1217.mps}
\end{subfigure}
\end{subfigure}
\caption{Group $\PGL_2(11)$: orbit 12.13 of size 2.}
\label{fig:PGL_2_11}
\end{figure}

$$
\beta = \frac{P_3^3}{cR} = \frac{Q_2^2\cdot Q_1}{cR},
$$
where
\begin{xflalign*}
c = {}
  & 2^8\cdot 3^3\cdot \omega^7\cdot (1-2\omega)^5\cdot (1-3\omega)^{11}
  \qquad\qquad\qquad\qquad
  \qquad\qquad
  \omega=\frac{1\pm\sqrt{5}}{2} \\
P_3 = {}
  & z^4 - 2 (1-2\omega)^3 z^3 - 2^2\omega (1-2\omega)^3 (3+\omega)^2 z^2 + {} \\
  & {} + 2\cdot 3^2\cdot 13 \omega^3 (1-2\omega)^4 z - 3 \omega^5 (1-2\omega)^4 (7-\omega) (5-8\omega) (11+\omega) \\
Q_2 = {}
  & z^5 - 7 (1-2\omega) z^4 - 2\omega^3 (1-2\omega) (127-12\omega) z^3 + {} \\
  & {} + 2\cdot 3^3\cdot 23 \omega^3 (1-2\omega)^2 z^2 - 3^2\omega^5 (1-2\omega)^2 (360-1469\omega) z - {}\\
  & {} - 3^3\cdot 13 \omega^6 (1-2\omega)^9\\
Q_1 = {}
  & z^2 - 2^4 (1-2\omega) z - \omega (1-2\omega) (221-19\omega)\\
R = {}
  & z^2 - \omega^3 (1-2\omega)^7.
\end{xflalign*}

%%%%%%%%%%%%%%%%%%%%%%%%%%%%%%%%%%%%%%%%%%

\subsectionsep{$n=13$, $(3^4 1^1 \,|\, 2^4 1^5 \,|\, 13^1)$, $ER=\PSL_3(3)$}

An ordinary tree (orbit 13.1, see~\cite{AdrZvo-2013}), the orbit consists of four dessins, see fig.~\ref{fig:PSL_3_3}.
In \cite{CasCou-1999} a one-parameter family of {\it Zolotarev polynomials}\footnote{Zolotarev polynomial is a complex polynomial with three finite critical values. Zolotarev polynomials relate to Shabat polynomials in the same way as Fried pairs (covers of $\PP^1$ ramified over four points) relate to Belyi pairs.} is calculated. The Shabat polynomial of this plane tree can be obtained from the family of Zolotarev polynomials by specializing its parameter $T$.

\begin{figure}[!ht]
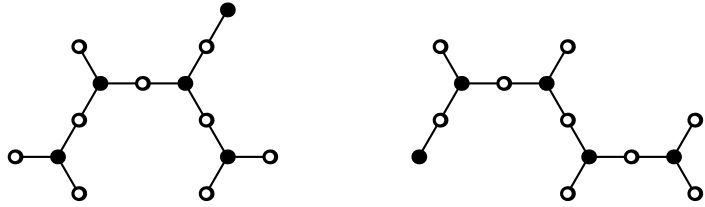

\centering
{\includegraphics[scale=0.8]{images/wt-1301.mps}
\hspace{1.5cm}
\includegraphics[scale=0.8]{images/wt-1302.mps}}
\caption{Group $\PSL_3(3)$: orbit 13.1 of size 4.}
\label{fig:PSL_3_3}
\end{figure}

The field of definition of the Shabat polynomials turns out to be the same as the field of definition of the family of Zolotarev polynomials. it can be defined as $K=\Q(a)$, where $a^4+13a^2+13=0$.
Note that $\eps = (a^2+2)/3$ satisfies $\eps^2+3\eps-1=0$, i.e. it is a unity of the field.

The Galois group of the field $K$ is cyclic of order 4, and hence it contains all the roots of the equation $a^4+13a^2+13=0$; the roots are $\{ a, a/\eps, -a, -a/\eps \}$.

To make our formulas more compact let us consider the equation
$\omega^4+\omega^3+2\omega^2-4\omega+3=0$; its roots are
$$
\begin{array}{rcl}
\omega_1 &=& (a^2 - 3a + 5)/6 \\
\omega_2 &=& (a^3 - a^2 + 11a - 8)/6 \\
\omega_3 &=& (a^2 + 3a + 5)/6 \\
\omega_4 &=& (-a^3 - a^2 - 11a - 8)/6 \\
\end{array}
$$
Then $a=\omega_3-\omega_1$, $\eps=\omega_1+\omega_3-1$ and the Belyi function can be written as
$$
\beta = \frac{P_3^3 z}{c} = \frac{Q_2^2\cdot Q_1}{c} + 1,
$$
where
$$
c = -2^6\eps^{11}\omega_1^{15}\omega_2^3\omega_3^{15}\omega_4^3 =
    -2^6\cdot 3^3\eps^{11}\omega_1^{12}\omega_3^{12},
$$
\begin{xflalign*}
\qquad
P_3 ={}& \left(z + a\omega_1^3\right)\left(z^3
  + a\omega_1\omega_2\omega_3(2 - \omega_1 + \omega_2 - \omega_3)z^2 +{}\right. & \\
&
  \qquad\qquad
  \left.{} + a^2\eps\omega_1\omega_3(2 - 3\omega_1 + 3\omega_2 + 3\omega_3)z
  + 2^3a^3\eps^3\right),
\end{xflalign*}
\begin{xflalign*}
\qquad
Q_2 ={}& z^4 + 2\cdot 5\eps\omega_1 z^3 +
  7\eps^2\omega_1\omega_3\omega_4(2-2\omega_1+\omega_3) z^2 +{} & \\
&
  \qquad{} + 2^2\eps^3(1 - \omega_1)(1 + 12\omega_1 + 6\omega_2 - 12\omega_3)z +
  2^3\eps^5\omega_1^3,
\end{xflalign*}
\begin{xflalign*}
\qquad
Q_1 ={}& \left( z - \omega_1^6/\eps \right)
  \left(z^4 - 2\omega_1\omega_2\omega_3(4 + 5\omega_1 + 3\omega_2 - 4\omega_3) z^3 -{}
  \right. & \\
&
  \qquad {}- \eps^2\omega_1\omega_3(2-\omega_1+\omega_2+\omega_3)
  (16+21\omega_1+14\omega_2+7\omega_3) z^2 -{}\\
&
  \qquad {}- 2\eps^5(6+8\omega_1+2\omega_2+7\omega_3)
  (16+4\omega_1+17\omega_2+5\omega_3) z -{}\\
&
  \qquad \left.{}- \eps^2\omega_1^3\omega_2^3\omega_3^{15}\omega_4^3\right).
\end{xflalign*}

%%%%%%%%%%%%%%%%%%%%%%%%%%%%%%%%%%%%%%%%%%

\subsectionsep{$n=14$, $(3^4 1^2 \,|\, 2^6 1^2 \,|\, 13^1 1^1)$, $ER=\PSL_2(13)$}

A weighted tree (orbit 14.1, see~\cite{AdrZvo-2013}), the orbit consists of one dessin, see fig.~\ref{fig:PSL_2_13}.
The Belyi function is the covering map $X_0(13)\to X(1)$ of the classical modular curves; it was calculated by F.\,Klein~\cite{Klein-1878}.

\begin{figure}[!ht]
\centering
\begin{subfigure}[b]{0.45\textwidth}
\centering
\includegraphics[angle=90,scale=0.8]{images/wt-1401.mps}
\vspace{0.3cm}
\end{subfigure}
\begin{subfigure}[b]{0.45\textwidth}
\centering
\includegraphics[width=6cm]{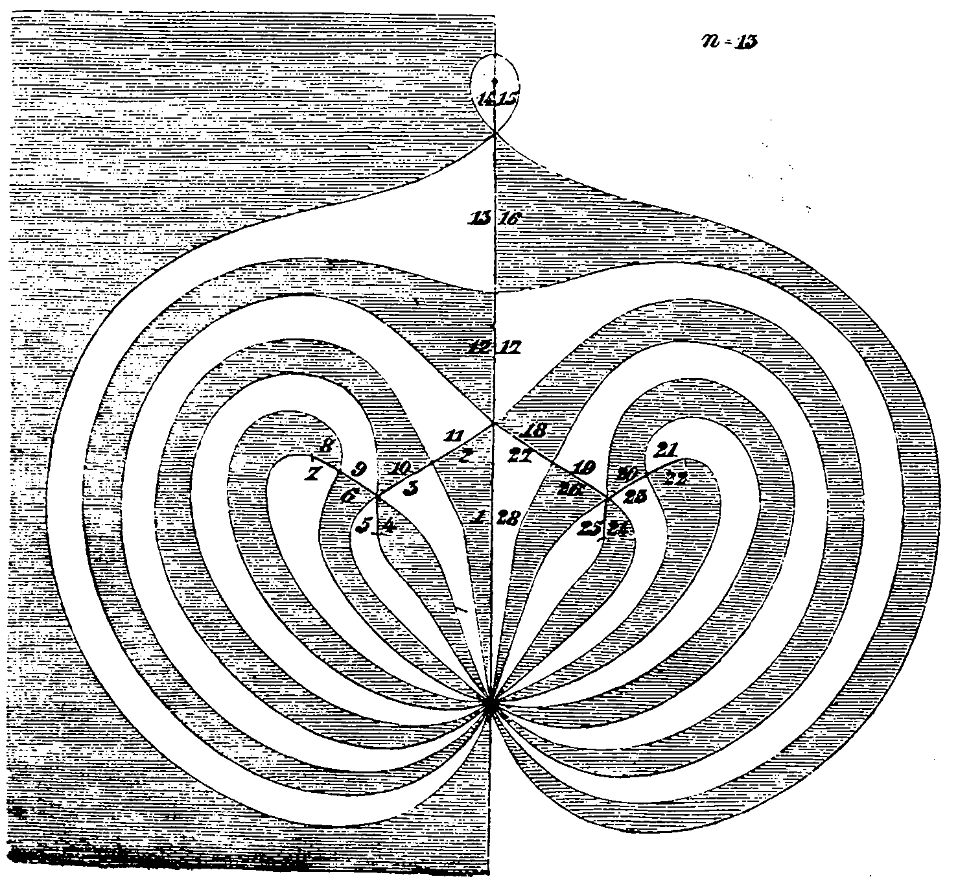}
\end{subfigure}
\caption{Group $\PSL_2(13)$: orbit 14.1 of size 1 and a drawing by F.\,Klein.}
\label{fig:PSL_2_13}
\end{figure}

$$
\beta = \frac{(z^4+7z^3+20z^2+19z+1)^3(z^2+5z+13)}{1728z}
$$
$$
\beta - 1 = \frac{(z^6+10z^5+46z^4+108z^3+122z^2+38z-1)^2(z^2+6z+13)}{1728z}.
$$

%%%%%%%%%%%%%%%%%%%%%%%%%%%%%%%%%%%%%%%%%%

\subsectionsep{$n=14$, $(3^4 1^2 \,|\, 2^7 \,|\, 12^1 1^2)$, $ER=\PGL_2(13)$}

A weighted tree (orbit 14.2, see~\cite{AdrZvo-2013}), the orbit consists of two dessins, see fig.~\ref{fig:PGL_2_13}.

\begin{figure}[!ht]
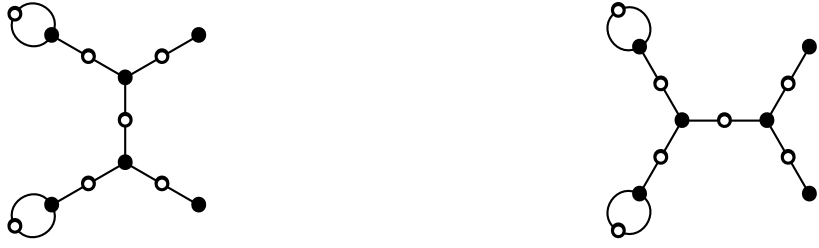

\centering
\begin{subfigure}[b]{0.45\textwidth}
\centering
\includegraphics[angle=90,scale=0.8]{images/wt-1402.mps}
\end{subfigure}
\begin{subfigure}[b]{0.45\textwidth}
\centering
\includegraphics[angle=90,scale=0.8]{images/wt-1403.mps}
\end{subfigure}
\caption{Group $\PGL_2(13)$: orbit 14.2 of size 2.}
\label{fig:PGL_2_13}
\end{figure}

The Belyi function is defined over $\Q(\sqrt{3})$; the number $\varepsilon=2\pm\sqrt{3}$ is a unit of the field.
Функция Белого оказывается определена над $\Q(\sqrt{3})$; число $\varepsilon=2\pm\sqrt{3}$ является единицей поля.
$$
\beta = \frac{P_3^3\cdot P_1}{cR} = \frac{Q_2^2}{cR},
$$
where
\begin{xflalign*}
c ={}&\frac{\omega^{16} (1-\omega)^{9}}{\varepsilon^4 (5-\omega)}
\qquad\qquad
\varepsilon=1+\omega
\qquad\qquad
\omega = 1\pm\sqrt{3}
\end{xflalign*}
\begin{xflalign*}
P_3 = {}
& z^4 + \omega^2 (1-\omega) (5-\omega) z^3 + \omega^2 (1-\omega)^3 (3-4\omega) z^2 + {}\\
& {} + \omega^2 (1-\omega)^3 (1-2\omega) (1-3\omega) z + \varepsilon^2 (1-\omega)^3 (35-19\omega)
\end{xflalign*}
\begin{xflalign*}
P_1 = {} & z^2 + (1-\omega)^2 = z^2 + 3
\end{xflalign*}
\begin{xflalign*}
Q_2 ={}
& z^7 + \varepsilon (1-\omega)^3 (5-\omega) z^6 + \varepsilon (1-\omega)^2 (51+7\omega) z^5 + {}\\
&{} + 5 \varepsilon^2 (1-\omega)^3 (47-10\omega) z^4 + \varepsilon^2 (1-\omega)^4 (301-12\omega) z^3 + {}\\
& {} + \varepsilon^2 (1-\omega)^5 (261+34\omega) z^2 + \varepsilon^2 (1-\omega)^6 (157+28\omega) z + {}\\
& {} + \varepsilon^3 (1-\omega)^7 (115-27\omega)
\end{xflalign*}
\begin{xflalign*}
R ={}& 13z^2 + \omega^2 (1-\omega)^5 (5-\omega) z + \varepsilon (1-\omega)^2 (5-\omega) (31+10\omega).
\end{xflalign*}

%%%%%%%%%%%%%%%%%%%%%%%%%%%%%%%%%%%%%%%%%%

\subsectionsep{$n=17$, $(3^5 1^2 \,|\, 2^8 1^1 \,|\, 15^1 1^2)$, $ER=\PSL_2(16)$}

A weighted tree (orbit 17.1, see~\cite{AdrZvo-2013}), the orbit consists of one dessin, see fig.~\ref{fig:PSL_2_16}.

\begin{figure}[!ht]
\centering
\includegraphics[angle=-90,scale=0.8]{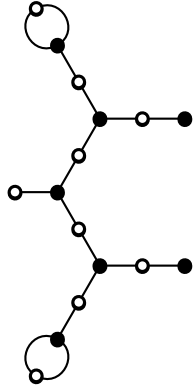}
\caption{Group $\PSL_2(16)$: orbit 17.1 of size 1.}
\label{fig:PSL_2_16}
\end{figure}

$$
\beta = \frac{(z^5 + 6z^4 + 48z^3 + 144z^2 + 432z + 288)^3 (z^2 + 6z + 24)}{2^{14}\cdot 3^6 (2z^2 + 3z + 48)}
$$
$$
\beta - 1 = \frac{z Q_2^2}{2^{14}\cdot 3^6 (2z^2 + 3z + 48)} + 1
$$
$$
Q_2 = z^8 + 12z^7 + 120z^6 + 720z^5 + 3600z^4 + 12096z^3 + 32832z^2 + 51840z + 51840.
$$

%%%%%%%%%%%%%%%%%%%%%%%%%%%%%%%%%%%%%%%%%%

\subsectionsep{$n=20$, $(3^6 1^2 \,|\, 2^9 1^2 \,|\, 18^1 1^2)$, $ER=\PGL_2(19)$}

A weighted tree (orbit 20.1, see~\cite{AdrZvo-2013}), the orbit consists of three dessins, see fig.~\ref{fig:PGL_2_19}.

\begin{figure}[!ht]
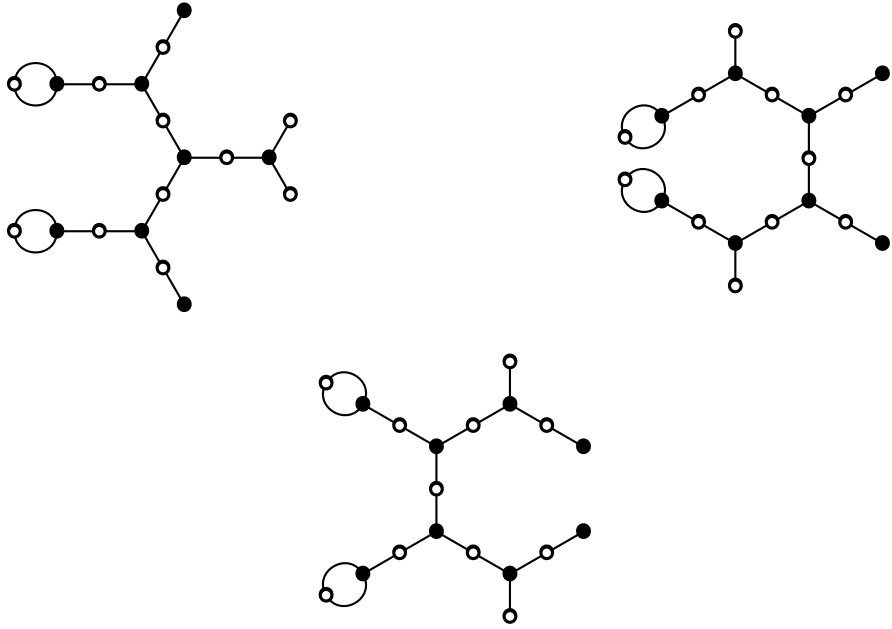

\centering
\begin{subfigure}[b]{0.45\textwidth}
\centering
\includegraphics[angle=90,scale=0.8]{images/wt-2001.mps}
\end{subfigure}
\begin{subfigure}[b]{0.45\textwidth}
\centering
\includegraphics[angle=90,scale=0.8]{images/wt-2002.mps}
\vspace{0.25cm}
\end{subfigure}
\vskip 0.5cm
\includegraphics[angle=90,scale=0.8]{images/wt-2003.mps}
\caption{Group $\PGL_2(19)$: orbit 20.1 of size 3.}
\label{fig:PGL_2_19}
\end{figure}

The Belyi function is defined over $K=\Q(a)$, where $a^3 - 3a + 1$.
The numbers $\varepsilon_1=a$ and $\varepsilon_2=a-1$ are units of the field.
The Galois group of the field $K$ is cyclic of order 3.

$$
\beta = \frac{P_3^3 P_1}{cR} = \frac{Q_2^2 Q_1}{cR} + 1,
$$
where
$$
c = -2^{16} \varepsilon_1^{50} \varepsilon_2^{-24} (a+1)^9 (a-3)^{19} \\
$$
$$
R = z^2 + \varepsilon_1 \varepsilon_2^6 (a+1)^{11} \\
$$
\begin{xflalign*}
P_3 = {} & z^6 - 2 \varepsilon_1^2 (a+1)^2 (a^2-a+3) z^5 +{} \\
  & {} + \varepsilon_1^2 \varepsilon_2^2 (a+1)^2 (3a-2) (3a+5)^2 (7a-3) z^4 - {} \\
  & {} - 4 \varepsilon_1^6 (a+1)^3 (30a^2+167a+22) z^3 - {} \\
  & {} - \varepsilon_1^{10} \varepsilon_2^{-3} (a+1)^4 (3a-4) (a^2+2a-9) (34a^2+8a+125) z^2 - {} \\
  & {} - 2 \varepsilon_1^{10} (a+1)^5 (1540a^2+5229a-206) z + {} \\
  & {} + \varepsilon_1^{14} \varepsilon_2^{-5} (a+1)^5 (3a-5) (11a-8) (a^2-6a+15) (5a^2-48a-21)
\end{xflalign*}
\begin{xflalign*}
P_1 = {} & z^2 - 2 \varepsilon_1 \varepsilon_2^2 (a+1)^5 z + \varepsilon_1^3 \varepsilon_2^2 (a+1)^5 (2a+1) (4a+3)
\qquad\qquad\qquad\quad
\end{xflalign*}
\begin{xflalign*}
Q_2 = {} & z^9 - \varepsilon_1^2 (a+1) (14a+9) z^8 + {} \\
  {} + {} & 4 \varepsilon_1^6 \varepsilon_2^{-6} (a+1) (3a-5) (5a-8) (18a-5) z^7 + {} \\
  {} + {} & 4 \varepsilon_1^5 \varepsilon_2^{-1} (a+1)^3 (12a^2-975a+397) z^6 + {} \\
  {} + {} & 2 \varepsilon_1^9 \varepsilon_2^3 (a+1)^4 (a^2-9a+6) (2a^2-4706a-8847) z^5 + {} \\
  {} + {} & 2 \varepsilon_1^{14} \varepsilon_2^{-5} (a+1)^4 (6a+5) (2a^2-11a-12) (11a^2+1116a-1729) z^4 - {} \\
  {} - {} & 4 \varepsilon_1^{15} \varepsilon_2^{-6} (a+1)^3 (3a-2) (6348a^2-115767a+151928) z^3 - {} \\
  {} - {} & 4 \varepsilon_1^{14} \varepsilon_2^{-3} (a+1)^4 (19a-7) (2a^2+15a+33) (25a^2-991a+3062) z^2 - {} \\
  {} - {} & \varepsilon_1^{19} (a+1)^4 (11a+21) (4a^2-28a-49) (1249a^2+4040a+23231) z + {} \\
  {} + {} & \varepsilon_1^{19} \varepsilon_2^2 (a+1)^{15} (a^2-3a-8) (a^2+3a+4) (7a+5) \cdot {} \\
  & \qquad\qquad\qquad\qquad\qquad\qquad\qquad {} \cdot (a^2-26a+40)(7a^2-31a+56)
\end{xflalign*}
\begin{xflalign*}
Q_1 = {} & z^2 - 2^3 \varepsilon_1^3 (a+1) z + \varepsilon_1^5 (a+1) (5a^2+40a+84).
\qquad\qquad\qquad\qquad\quad
\end{xflalign*}

%%%%%%%%%%%%%%%%%%%%%%%%%%%%%%%%%%%%%%%%%%

\subsectionsep{$n=24$, $(3^6 1^6 \,|\, 2^{12} \,|\, 23^1 1^1)$, $ER=\M_{24}$}

A weighted tree (orbit 24.1, see~\cite{AdrZvo-2013}), the orbit consists of two dessins, see fig.~\ref{fig:M_24a}. The Belyi function is calculated in \cite{HoyMat-1986}. The dessin is also used in \cite{Mon-2014} to demonstrate the method of calculation based on the modular functions techniques.

\begin{figure}[!ht]
\centering
\includegraphics[scale=0.8]{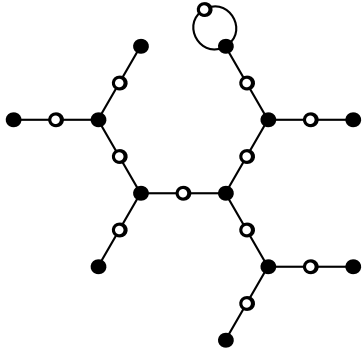}
\caption{Group $\M_{24}$: orbit 24.1 of size 2.}
\label{fig:M_24a}
\end{figure}

$$
\beta = \frac{P_3^3\cdot P_1}{cz} = \frac{Q_2^2}{cz},
$$
where
\begin{xflalign*}
c={}&-2^{38}\cdot 3^3(1-\omega)^9(1+\omega)^4&
\omega=\frac{1\pm\sqrt{-23}}{2}
\end{xflalign*}
\begin{xflalign*}
P_3={}&
  z^6 - 2\cdot 3\cdot 7 z^5 -(1-\omega)^3(23+16\omega)z^4 -{}&\\
&
  {}-2(1+\omega)(616-423\omega)z^3 - 
  (1-\omega)(29-6\omega)(37-139\omega)z^2 +{}\\
&
  {}+2(1-\omega)^3(1099+45\omega)z + (1-\omega)(1+\omega)^5(1+2\omega)
\end{xflalign*}
\begin{xflalign*}
P_1={}&
  z^6 - 2\cdot 29 z^5 - 5(1-3\omega)(9-17\omega)z^4 +{}&\\
&
  {}+2\omega(1+\omega)(1043+348\omega)z^3+
  (1+\omega)(23675-23163\omega)z^2 -{}\\
&
  {}-(1+\omega)(20-13\omega)(37+15\omega)(67-47\omega)z +{}\\
&
  {}+23^2(1-\omega)(1+\omega)^5(1+2\omega)
\end{xflalign*}
\begin{xflalign*}
Q_2={}&
  z^{12} - 2^2\cdot 23 z^{11} + 2\cdot 23(1-\omega)(1+14\omega)z^{10} -{}&\\
&
  {}-2^2\cdot 23(1+\omega)(245-141\omega)z^9 +{}\\
&
  {}+23(1+\omega)(3-\omega)(7427-1236\omega)z^8 -{}\\
&
  {}-2^3\cdot 3^2\cdot 23(1-\omega)^2(1+7\omega)(25+43\omega)z^7 +{}\\
&
  {}+23(1+\omega)^2(4+3\omega)(4-21\omega)(211-779\omega)z^6 +{}\\
&
  {}+2^3\cdot 23(1+\omega)^2(470573+119123\omega)z^5 +{}\\
&
  {}+23(1-\omega)^2(1+\omega)^2(13+6\omega)(59-74\omega)(587+162\omega)z^4 +{}\\
&
  {}+2\cdot 3\cdot 23(1-\omega)^3(1+\omega)^2(7-2\omega)(7420-25173\omega)z^3 +{}\\
&
  {}+2\cdot 23(1-\omega)^5(1+\omega)(5-\omega)(11641-65466\omega)z^2 -{}\\
&
  {}-23(1+\omega)^7(2+\omega)(1+2\omega)^2(1099+45\omega)z -{}\\
&
  {}-23(1-\omega)^2(1+\omega)^{10}(1+2\omega)^2.
\end{xflalign*}

%%%%%%%%%%%%%%%%%%%%%%%%%%%%%%%%%%%%%%%%%%

\subsectionsep{$n=24$, $(3^8 \,|\, 2^8 1^8 \,|\, 23^1 1^1)$, $ER=\M_{24}$}

A weighted tree (orbit 24.2, see~\cite{AdrZvo-2013}), the orbit consists of two dessins, see fig.~\ref{fig:M_24b}.

\begin{figure}[!ht]
\centering
\includegraphics[scale=0.8]{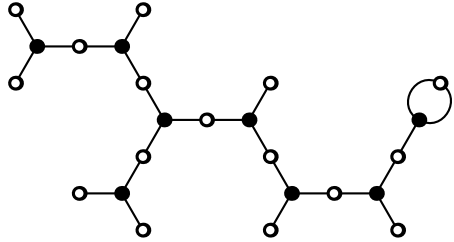}
\caption{Group $\M_{24}$: orbit 24.2 of size 2.}
\label{fig:M_24b}
\end{figure}

%Функцию Белого можно записать в виде
$$
\beta = \frac{P_3^3}{cz} = \frac{Q_2^2\cdot Q_1}{cz}
$$
\begin{xflalign*}
c={}&
-2^{18}\cdot 3^3\cdot \omega^{12}\cdot (1+\omega)^{10}&
\omega=\frac{1\pm\sqrt{-23}}{2}
\end{xflalign*}
\begin{xflalign*}
P_3={}&
  z^8 - 23(1+\omega)  z^7
  + 2^2\cdot 23\omega(5+2\omega)  z^6 +{}&\\
&
  {}+ 2^3\cdot 23(8+\omega)(15-2\omega)  z^5
  - 2\cdot 23\omega(2-\omega)(380+149\omega)  z^4 -{}\\
&
  {}- 2^3\cdot 23(2-\omega)(2+5\omega)(63-44\omega)  z^3 -{}\\
&
  {}- 2^2\cdot 23(2-\omega)(14-\omega)(14-169\omega)  z^2 -{}\\
&
  {}- 2^2\cdot 23(2-\omega)^3(61-197\omega)  z
  - 23(2-\omega)^4(3-2\omega)^2
\end{xflalign*}
\begin{xflalign*}
Q_2={}&
  z^8 - 2^2\cdot 5(1+\omega)  z^7
  + 2^2\cdot 17\omega(5+2\omega)  z^6 -{}&\\
&
  {}- 2^3\cdot 3\cdot 7(2-3\omega)(1-5\omega)  z^5
  - 2\cdot 11\omega(2-\omega)(380+149\omega)  z^4 -{}\\
&
  {}- 2^6(2-\omega)(2+5\omega)(63-44\omega)  z^3 -{}\\
&
  {}- 2^2\cdot 5(2-\omega)(14-\omega)(14-169\omega)  z^2 -{}\\
&
  {}- 2^3(2-\omega)^3(61-197\omega)  z
  +(2-\omega)^4(3-2\omega)^2
\end{xflalign*}
\begin{xflalign*}
Q_1={}&
  z^8 - 29(1+\omega) z^7
  -(3-\omega)(833-114\omega) z^6 +{}&\\
&
  {}+ 5(1+\omega)(3-4\omega)(335+64\omega)  z^5 -{}\\
&
  {}-(2-\omega)(4-3\omega)(3-7\omega)(593-696\omega)  z^4 -{}\\
&
  {}- 2^2\omega(8-3\omega)(9629-17877\omega)  z^3 -{}\\
&
  {}-(2-\omega)^2(7-3\omega)(41+83\omega)(524-387\omega) z^2 +{}\\
&
  {}+ 2^2(2-\omega)(8-\omega)(332225-49341\omega)  z
  - 23^3\omega^6(2-\omega)^2.
\end{xflalign*}

%%%%%%%%%%%%%%%%%%%%%%%%%%%%%%%%%%%%%%%%%%

\subsectionsep{$n=32$, $(3^{10} 1^2 \,|\, 2^{12} 1^8 \,|\, 31^1 1^1)$, $ER=\ASL_5(2)$}

A weighted tree (orbit 32.1, see~\cite{AdrZvo-2013}), the orbit consists of 6 dessins, see fig.~\ref{fig:ASL_5_2}.

\begin{figure}[!ht]
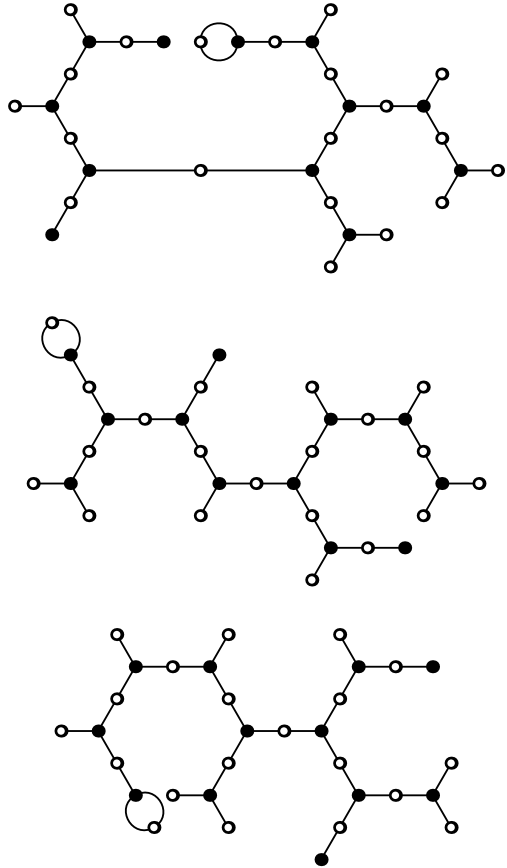

\centering
\includegraphics[scale=0.7]{images/wt-3201.mps}
\vskip 0.5cm
\includegraphics[scale=0.7]{images/wt-3202.mps}
\vskip 0.5cm
\includegraphics[scale=0.7]{images/wt-3203.mps}
\caption{Group $\ASL_5(2)$: orbit 32.1 of size 6.}
\label{fig:ASL_5_2}
\end{figure}

One can use GAP\footnote{https://www.gap-system.org/} to check that there exist 6 conjugacy classes of elements of order 31 in $\ASL_5(2)$. For any two of then $C_1$, $C_2$ there exists $k\in (\Z/31\Z)^*$ such that any element of $C_1$ is the $k$-th power of an element of $C_2$. There are 6 irrational numbers among the values of irreducible characters on these classes; they are of the form $e^k_{31}+e_{31}^{2k}+e_{31}^{4k}+e_{31}^{8k}+e_{31}^{16k}$, where $e_{31}$ is a primitive root of unity of degree 31. The numbers are the roots of the polynomial $a^6+a^5+3a^4+11a^3+44a^2+36a+32$. The Galois group of the polynomial is cyclic of order 6.

It follows from the general theory that the Belyi functions of these dessins are defined over the splitting field of $a^6+a^5+3a^4+11a^3+44a^2+36a+32$; our calculations confirm it. This field of definition was also calculated in 2016 by J.\,Voight.
%\footnote{The authors are indebted to A.\,Zvonkin for this information.}.
The formulas of the Belyi function are too unwieldy to bring them here; there reader can find them in the supplementary materials.

\bigskip

{\bf Supplementary materials.} 
The Belyi functions presented in this paper are available as SageMath files in the GitHub repository \url{https://github.com/nadrianov/supplementary}.

\end{document}